\documentclass[a4paper,12pt]{amsart}

\usepackage{amssymb, amscd}
\usepackage{latexsym,epsfig}
\usepackage[all]{xy}
\usepackage{pst-all}
\numberwithin{equation}{section}
\def\today{\ifcase\month\or Jan\or Febr\or  Mar\or  Apr\or May\or Jun\or  Jul\or Aug\or  Sep\or  Oct\or Nov\or  Dec\or\fi \space\number\day, \number\year}

\newcommand{\CC}{\mathbb C}
\newcommand{\EE}{\mathbb E}

\newcommand{\PP}{\mathbb P}
\newcommand{\QQ}{\mathbb Q}

\newcommand{\ZZ}{\mathbb Z}
\newcommand\M[1]{{\mathcal M}_{#1}}

\newcommand\A[1]{{\mathcal A}_{#1}}
\newcommand\tA[1]{\widetilde{\mathcal A}_{#1}}

\newcommand\SL[2]{{\mathrm{SL}}({#1},{#2})}
\newcommand\GL[2]{{\mathrm{GL}({#1},{#2})}}
\newcommand\PGL[2]{{\mathrm{PGL}({#1},{#2})}}
\newcommand\Sp[2]{{\mathrm{Sp}}({#1},{#2})}
\newcommand{\Sym}{{\mathrm{Sym}}}

\numberwithin{equation}{section}
\newtheorem{theorem}{Theorem}[section]
\newtheorem{lemma}[theorem]{Lemma}

\newtheorem{conclusion}[theorem]{Conclusion}
\newtheorem{corollary}[theorem]{Corollary}

\newtheorem{definition-lemma}[theorem]{Definition-Lemma}

\theoremstyle{definition}

\newtheorem{example}[theorem]{Example}

\theoremstyle{remark}
\newtheorem{remark}[theorem]{Remark}

\begin{document}

\title[Covariants of binary sextics and Siegel modular forms]{Covariants of Binary Sextics and 
Vector-valued Siegel Modular Forms of Genus Two}

\author{Fabien Cl\'ery}
\address{ Max-Planck-Institut f\"ur Mathematik,
Vivatsgasse 7,
53111 Bonn,
Germany}
\email{cleryfabien@gmail.com}

\author{Carel Faber}
\address{Mathematisch Instituut, Universiteit Utrecht,
P.O. Box 80010,
3508 TA Utrecht,
The Netherlands}
\email{C.F.Faber@uu.nl}

\author{Gerard van der Geer}
\address{Korteweg-de Vries Instituut, Universiteit van
Amsterdam, Postbus 94248,
1090 GE  Amsterdam, The Netherlands.}
\email{geer@science.uva.nl}

\begin{abstract}
We extend Igusa's description of the relation between invariants of binary
sextics and Siegel modular forms of degree $2$ 
to a relation between covariants and vector-valued Siegel modular forms
of degree $2$. We show how this relation can be used to effectively 
calculate the Fourier expansions of Siegel modular forms of degree $2$.
\end{abstract}

\maketitle
%%%%%%%%%%%%%%%%%%%%%%%%%%%%%%%%%%%%%%%%
%\centerline{\today}
%%%%%%
\begin{section}{Introduction}
In his 1960 papers Igusa explained the relation between the invariant theory
of binary sextics and scalar-valued Siegel modular forms of degree (or genus) $2$.
The relation stems from the fact that the moduli space
$\M{2}$ of curves of genus $2$ admits two descriptions. 
The moduli space $\M{2}$
has a classical description in terms of the invariant theory of binary sextics.
But via the Torelli morphism $\M{2}$ 
is an open part of a Shimura variety, 
the moduli space $\A{2}$ of principally polarized abelian surfaces.
Therefore we have two descriptions of natural vector 
bundles on our moduli space and thus two descriptions of their sections. 
The purpose of this paper is to extend the correspondence
given by Igusa between invariants of binary sextics and
scalar-valued Siegel modular forms of degree $2$ to a correspondence between covariants
and vector-valued Siegel modular forms of degree~$2$. We
give a description of vector-valued Siegel modular forms of 
degree $2$ in terms of covariants
of the action of ${\rm SL}(2,{\CC})$ on the space of binary sextics. 
%Grace and Young showed that the ring of covariants of binary sextics
%is finitely generated and generated by $26$ elements.

Since the complement of the Torelli image of ${\mathcal M}_2$ in 
${\mathcal A}_2$ is the space ${\mathcal A}_{1,1}$
of decomposable abelian surfaces (products of elliptic curves) 
we have to analyze the orders of modular forms described by 
covariants along this locus. The scalar-valued modular
forms $\chi_{10}$ and its square root $\chi_5$ 
that vanish with multiplicity $2$ and $1$ along this locus play
a key role.

We shall see that the `first' vector-valued modular form that appears
on $\Sp{4}{\ZZ}$, the form 
$\chi_{6,8}$ of weight $(6,8)$, 
corresponds in some sense to the `universal' binary sextic.
The relation between covariants and modular forms allows us to use 
constructions from representation theory to construct vector-valued 
Siegel modular forms. 
In fact, we shall show that up to multiplication by a suitable 
power of $\chi_{10}$, every vector-valued Siegel modular form 
of degree $2$ can be obtained by applying a
representation-theoretic construction to the form $\chi_{6,8}$. 
In this sense, our form $\chi_{6,8}$
may be called the `universal vector-valued Siegel modular form of degree $2$'.
In fact, in practice we work with the form 
$\chi_{6,3}=\chi_{6,8}/\chi_{5}$
which is a modular form with a character.

We show that our constructions can be used to efficiently calculate the 
Fourier series of vector-valued Siegel modular forms. 
We illustrate this with a number of significant examples
where we use these Fourier series to calculate eigenvalues of Hecke operators and check instances of
Harder type congruences.

Similar ideas will be worked out for the case of genus $3$ in a forthcoming paper.

\noindent
{\sl Acknowledgement.} The first author thanks the Max-Planck-Institute
for the hospitality he is enjoying there.
\end{section}
%%%%%%%%%%%%%%%%%%%%%%%%%%%%%%%%%%%%%%%%
\begin{section}{Siegel Modular Forms}
The Siegel modular group $\Gamma_g=\Sp{2g}{\ZZ}$ of degree (or genus) $g$ 
acts on the Siegel upper half space 
$$
\mathfrak{H}_g=\{ \tau \in {\rm Mat}(g \times g, {\CC}):
\tau^t=\tau, {\rm Im}(\tau) > 0\}
$$ 
of degree $g$ by
$$
\tau \mapsto \gamma\cdot \tau =(a\tau+b)(c\tau+d)^{-1} \qquad \text{\rm for 
$\gamma=\left(\begin{matrix} a & b \\ c & d \\ \end{matrix} \right)\in \Gamma_g\, .$}
$$
A scalar-valued Siegel modular form of degree $g>1$ and 
weight $k$ is a holomorphic function
$f \colon \mathfrak{H}_g \to {\CC}$ satisfying 
$$
f(\gamma \cdot \tau)= \det(c\tau+d)^k f(\tau)
$$
for all $\gamma \in \Gamma_g$. (For $g=1$ we also need a growth condition at infinity.) 
If $W$ is a finite-dimensional complex vector space
and $\rho: \GL{g}{\CC} \to {\rm GL}(W)$
 a representation, then a vector-valued 
Siegel modular form of degree $g>1$ and weight $\rho$
is a holomorphic map $f: \mathfrak{H}_g \to W$ such that 
$$
f(\gamma\cdot \tau)= \rho(c\tau+d) f(\tau)
$$
for all $\gamma \in \Gamma_g$.
Siegel modular forms can be interpreted as sections of vector bundles.
The moduli space $\A{g}({\CC})=\Gamma_g \backslash \mathfrak{H}_g$
of complex principally polarized abelian varieties carries an orbifold vector
bundle of rank $g$, the Hodge bundle ${\EE}$, defined by the factor of automorphy 
given by
$$
j(\gamma,\tau)= (c\tau+d)\, .
$$
Its determinant $L=\det({\EE})$ defines a line bundle. A scalar-valued 
Siegel modular form
of degree $g$ and weight $k$ can be seen as a section of $L^{\otimes k}$.
A vector-valued Siegel modular form of weight $\rho$ can be viewed as a 
section of the vector bundle ${\EE}_{\rho}$, that is defined by the factor of
automorphy $j(\gamma,\tau)=\rho(c\tau+d)$. The vector bundle ${\EE}$ and
the bundles ${\EE}_{\rho}$ extend in a canonical way to toroidal compactifications
of $\A{g}$ and the Koecher principle says that their sections do so too.

We are interested in the case $g=2$. An irreducible representation of $\GL{2}{\CC}$
is of the form $\Sym^j(V)\otimes \det(V)^{\otimes k}$ for 
$(j,k)\in {\ZZ}_{\geq 0}\times {\ZZ}$ with $V$ the standard 
representation of $\GL{2}{\CC}$. We denote the corresponding vector space of Siegel
modular forms by $M_{j,k}(\Gamma_2)$ or simply by $M_{j,k}$. It contains a
subspace of cusp forms denoted by $S_{j,k}(\Gamma_2)$ or simply by $S_{j,k}$. 
The scalar-valued modular forms
correspond to the case where $j=0$.

The scalar-valued modular forms of degree $2$ form a ring 
$R=\oplus_k M_{0,k}(\Gamma_2)$. Igusa showed in \cite{Igusa1962} that $R$ is generated by
modular forms
$E_4,E_6,\chi_{10}$, $\chi_{12}$ and $\chi_{35}$ of weight $4,6,10,12$ and $35$.
The subring $R^{\rm ev}$ of modular forms of even weight is the pure polynomial ring
generated by $E_4,E_6,\chi_{10}$ and $\chi_{12}$, and the form $\chi_{35}$ satisfies
a quadratic equation over $R^{\rm ev}$ expressing $\chi_{35}^2$ as a polynomial in
$E_4,E_6,\chi_{10}$ and $\chi_{12}$. 

The bi-graded $R$-module $M=\oplus_{j,k} M_{j,k}$ can also be made into a ring by using
the canonical projections $\Sym^{j_1}(V) \otimes \Sym^{j_2}(V) \to \Sym^{j_1+j_2}(V)$
defined by multiplying homogeneous polynomials of degree $j_1$ and $j_2$ in two
variables to define the multiplication. 

The ring $M$ of vector-valued modular forms is not finitely generated as shown 
by Grundh, see  \cite[p.\ 234]{vdG}. 

\end{section}
%%%%%%%%%%%%%%%%%%%%%%%%%%%%%%%%%%%%%%%%
\begin{section}{Invariants and Covariants of Binary Sextics}
%%%
Let $V$ be a $2$-dimensional vector space (over ${\CC}$), generated
by elements $x_1$ and $x_2$. We will denote by ${\rm Sym}^6(V)$ 
the space of binary sextics, where we write a binary sextic as
$$
f=\sum_{i=0}^6 a_i \, \binom{6}{i} \, x_1^{6-i}x_2^i\, . \eqno(1)
$$

The group $\GL{2}{\CC}$ acts on $V$
via $(x_1,x_2) \mapsto (ax_1+bx_2, cx_1+dx_2)$ for $(a,b;c,d) \in \GL{2}{{\CC}}$.
This induces an action of $\GL{2}{\CC}$ on ${\rm Sym}^6(V)$,
and similarly actions of $\PGL{2}{\CC}$ and $\SL{2}{\CC}$ 
on ${\mathcal X}={\PP}({\rm Sym}^6(V))$.
We take $(a_0:a_1: \ldots :a_6)$ as coordinates on the projectivized space
${\mathcal X}$.

The space ${\mathcal X}$ carries the natural ample line bundle 
${\mathcal L}={\mathcal O}_{\PP(\Sym^6(V))}(1)$. On ${\mathcal L}$ we have
an action of the group $\SL{2}{\CC}$ that is compatible with the action on
the projectivized space ${\mathcal X}$ of binary sextics, 
see \cite{Mumford-GIT}.
We can retrieve
${\mathcal X}$ as the projective scheme 
${\rm Proj} (\oplus_n {\rm H}^0({\mathcal X},{\mathcal L}^n))$. 

The ring  of invariants is defined as the graded ring
$$
I\colon=\oplus_{n \geq 0} {\rm H}^0({\mathcal X}, {\mathcal L}^n)^{\SL{2}{\CC}} \, .
$$ 
Thus an invariant can be seen as a polynomial in the coefficients $a_i$ of 
$f$ that is invariant under the action of $\SL{2}{\CC}$. 
The ring of these invariants was determined in the 19th century by work of Clebsch 
and Bolza  \cite{Clebsch, Bolza}.
It is generated by elements $A,B,C,D$ and $E$ of degrees $6, 12, 18, 30$ and
$45$ in the roots of the binary sextic. We normalize their degree by
taking the degree in the roots divided by $3$. 
The invariant $D$ is the discriminant
and $E$ satisfies a quadratic equation expressing $E^2$ in the even
invariants.
\bigskip

Given the representation of $\GL{2}{\CC}$ on $V$ and on $\Sym^6(V)$
we look at equivariant inclusions of  $\GL{2}{\CC}$-representations
$$
\iota_A: A \longrightarrow \Sym^d(\Sym^6(V))\, ,
$$
or equivalently at equivariant inclusions 
$$
\varphi: {\CC}  \longrightarrow \Sym^d(\Sym^6(V))\otimes A^{\vee}\, .
$$
The image  $\varphi(1)$ is denoted $\Phi=\Phi(\varphi)$ 
and called a covariant.

If $A=A[\lambda]$ is the irreducible representation of $\GL{2}{\CC}$ of
highest weight $\lambda=(\lambda_1,\lambda_2)$ with $\lambda_1 \geq \lambda_2$
then the covariant $\Phi$ can be viewed as a form of degree $d$ in the variables 
$a_0,\ldots,a_6$, of degree $\lambda_1-\lambda_2$ in $x_1$ and $x_2$
and may be denoted by $\Phi(d,\iota_A)$, or simply by $\Phi(d,A)$
if $\iota_A$ is unique.

\begin{example}\label{tautologicalf} The tautological $f$.
As a first example we consider the (tautological) case where 
$d=1$ and $A=\Sym^6(V)$ and we take $\iota_A={\rm id}_{\Sym^6(V)}$.
Then the covariant $\Phi=\varphi(1)$ can be viewed as the universal binary
sextic $f=\sum_J a_J x^J$ given by ~(1).
\end{example}
\begin{example}\label{example-d2}
As a second example we look at $d=2$. We have the isotypical decomposition
$$
\Sym^2(\Sym^6(V))=A[12,0]+A[10,2]+A[8,4]+A[6,6] \, .
$$
We can view $\Phi(2,A[12,0])$ as the square of the tautological $f$ given by (1),
while the covariant $\Phi(2,A[10,2])$ is the Hessian of the polynomial $f$ given by 
$f_{x_1x_1}f_{x_2x_2}-f_{x_1x_2}^2$ or 
equivalently by
$$
(a_0a_2-a_1^2)x_1^8+ 4(a_0a_3-a_1a_2)x_1^7x_2+ \cdots + (a_4a_6-a_5^2)x_2^8.
$$
The covariant $\Phi$ corresponding to $A[8,4]$ is given by
$$
\begin{aligned}
(a_0a_4-4a_1a_3+3a_2^2)x_1^4+(2a_0a_5-6a_1a_4+4a_2a_3)x_1^3x_2
+(a_0a_6-9a_2a_4+8a_3^2)x_1^2x_2^2 &\\
+(2a_1a_6-6a_2a_5+4a_3a_4)x_1x_2^3+(a_2a_6-4a_3a_5+3a_4^2)x_2^4& \, . \\
\end{aligned}
$$
 Moreover, $\Phi(2,A[6,6])$ is an invariant equal to 
$a_0a_6-6\, a_1a_5+15\, a_2a_4-10\, a_3^2$ and coincides
(up to a multiplicative scalar) with the invariant $A \in I$.

\end{example}
The covariants of the action of $\SL{2}{\CC}$ form a graded ring $C$. This ring
can be identified with the ring of invariants
$$
{\CC}[V_1 \oplus V_6]^{\SL{2}{\CC}}\, ,
$$
where we write $V_i=\Sym^i(V)$, see \cite[p.\ 55]{Springer}.
It was shown in the 19th century by Clebsch and others that the ring $C$
is generated by $26$ elements ($5$ invariants and $21$ other covariants) 
\cite[p.\ 296]{Clebsch}. Covariants corresponding to an irreducible representation
 $A[\lambda]$
with $\lambda_1+\lambda_2=6 \, d$ can be calculated by the so-called symbolic method.
We refer to \cite{G-Y}, \cite{Chipalkatti} and \cite[p.\ 214]{M-S}. 

\end{section}
%%%%%%%%%%%%%%%%%%%%%%%%%%%%%%%%%%%%%%%%
%%%%%%%%%%%%%%%%%%%%%%%%%%%%%%%%%%%%%%%%
\begin{section}{Covariants and Siegel Modular Forms}\label{covariantsandSMF}
The moduli space of curves of genus $2$ admits two different 
descriptions. First, using the Torelli morphism $t$ we can view $\M{2}$
as an open subspace of the moduli space $\A{2}$ of principally 
polarized abelian surfaces. The complement of the image $t(\M{2})$
in $\A{2}$ is the divisor $\A{1,1}$ of products of elliptic curves.
Over the complex numbers we can view $\M{2}({\CC})$ as an open suborbifold
of the orbifold $\A{2}({\CC})=\Gamma_2\backslash \mathfrak{H}_2$.

The second description of $\M{2}$ is as the GIT quotient associated
to the action of $\GL{2}{\CC}$
on the space $X=\Sym^6(V)$ of binary sextics.
A binary sextic $f(x_1,x_2)$ with non-vanishing discriminant 
determines a curve $C$ of genus $2$ with affine equation
$y^2=\tilde{f}(x)$ with $\tilde{f}(x)=f(x,1)$. Moreover, the curve $C$ comes
with a basis of the space of differentials given by $xdx/y$ and $dx/y$.

We let $\GL{2}{\CC}$ act on the pairs $(x,y)$ by
$$
x \mapsto (ax+b)/(cx+d), \qquad y \mapsto y/(cx+d)^3\, .
$$
This gives an action $dx \mapsto (ad-bc) \, dx /(cx+d)^2$ and hence
an action by $\det (A) \, A$ on our basis $xdx/y, dx/y$.

Using this we can identify the moduli space $\M{2}$ with the 
algebraic stack $[Y/\GL{2}{\CC}]$ with $Y$ the algebraic stack of curves with
a framed Hodge bundle: the objects of $Y$ are pairs $(\pi,\alpha)$ 
with  $\pi: C \to S$ a curve of genus $2$ and an isomorphism 
$\alpha: {\mathcal O}_S^{\oplus 2} {\buildrel \sim \over \longrightarrow} 
\pi_*\omega_{\pi}$.
The group $\GL{2}{\CC}$ acts via $(\pi,\alpha)\mapsto
A \cdot (\pi,\alpha)=(\pi, \det(A)\, A\circ \alpha)$. 
We thus have (see \cite{Vistoli} )
$$
q : [Y/\GL{2}{\CC}] {\buildrel \sim \over \longrightarrow } \M{2} \, .
$$ 

Let $X^0\subset X$ be the complement of the discriminant locus.
We have a natural identification of $Y$ with $X^0$ such that the above action
of $\GL{2}{\CC}$ on $Y$ corresponds to the action of $\GL{2}{\CC}$ on $X$
given by $A \circ f =  f(ax_1+bx_2,cx_1+dx_2)$. This implies the following corollary.

\begin{corollary} \label{pullbackE}
The pullback of the Hodge bundle ${\EE}$ on $\M{2}$ to $X^0\cong Y$ is the
$\GL{2}{\CC}$-equivariant bundle $V\otimes \det(V)$.
\end{corollary} 

This can be extended to include the locus of binary sextics with 
at least five distinct zeros. We write a binary sextic with five
zeros as $(x-\alpha)^2f_4$ with $f_4$ a binary quartic with non-zero
discriminant and $f_4(\alpha)\neq 0$. To this we associate the nodal
curve $C$ which is obtained by associating to $f_4$ the genus $1$ curve $C_1$
given by $w^2=\tilde{f}_4$ with $\tilde{f}_4=f_4(x,1)$ and identifying
the two points of $C_1$ lying over $x=\alpha$. It comes with two 
differential forms $dx/w$ and $dx/(x-\alpha)w$; the latter form has
poles with opposite residues in the two points of $C_1$ lying over $x=\alpha$.
The connection with the case of sextics with non-vanishing discriminant
is given by setting $y=w (x-\alpha)$. In this way we can associate
to a binary sextic with at least five zeros a nodal curve with a frame
of the Hodge bundle. The group $\GL{2}{\CC}$ acts by 
$$
x \mapsto  (ax+b)/(cx+d), \qquad w \mapsto w/(cx+d)^2 \, .
$$
The identification $[Y/\GL{2}{\CC}] \cong \M{2}$ can be extended
if we include for $Y$ the case of irreducible nodal curves of 
arithmetic genus $2$ with one node and replace $\M{2}$ by 
$\overline{\mathcal M}_2-\Delta_1$ with $\Delta_1$ the locus of reducible curves.
The conclusion of Corollary \ref{pullbackE} is still valid.
In a similar way it can  be extended to the case of binary sextics all
of whose zeros have multiplicity $\leq 2$. 
We denote by $X^s$  the open subset
of $X$ of stable sextics, that is, binary sextics 
none of whose zeros have multiplicity $\geq 3$. 

\smallskip

A scalar-valued Siegel modular form $F$ of weight $k$ on 
$\Gamma_2$ is by definition a
section of the line bundle $L^k$ with $L=\det({\EE})$ on $\tilde{\mathcal A}_2$,
the standard toroidal compactification (equal to $\overline{\mathcal M}_2$),
and its pullback to $X^s$ will give rise to an invariant 
section of $\det(V)^{3k}$ on
${\mathcal X}$, that is, an invariant $i_F$, hence an invariant on all of $X$. 
This gives a map from the ring of scalar-valued Siegel modular
forms $R$ to the ring of invariants
$$
i: R \longrightarrow I\, .
$$ Igusa defined this injective
map in a slightly different way in \cite{Igusa1967}.
By Igusa's map we can view $I$ as a ring of meromorphic Siegel modular forms.
The generators
$E_4,E_6,\chi_{10}$, $\chi_{12}$ and $\chi_{35}$ of $R$
correspond to  non-zero multiples of
$B$, $AB-3C$, $D$, $AD$ and $D^2E$ (see \cite[p.\ 848]{Igusa1967}).
We see that every invariant
defines a meromorphic modular form $F$ such that $\chi_{10}^m F$ is
holomorphic for an appropriate power $m$.
In other words we have inclusions
$$
R \subset I \subset R_{\chi_{10}} \, ,
$$
where $R_{\chi_{10}}$ is obtained by inverting $\chi_{10}$.
Moreover, the ideal of $R$ of cusp forms maps to $(D)$, the
ideal generated by the discriminant (cf.\ \cite[p.\ 845]{Igusa1967}).

The relation between invariants and scalar-valued Siegel modular forms
can be extended to a relation between covariants and 
vector-valued Siegel modular forms. A section $F$ of
${\rm Sym}^j({\EE})\otimes {\det}({\EE})^k$ 
will by pullback under $q$ give rise to a section  of 
${\rm Sym}^{j}(V) \otimes \det(V)^{j+3k}$, hence to a covariant of the action
of ${\rm SL}(2,{\CC})$ on ${\rm Sym}^6(V)$.
We get inclusions similar to those above 
$$
M \subset C \subset M_{\chi_{10}}
$$
with $M_{\chi_{10}}=M \otimes_R R_{\chi_{10}}$ 
fitting into a commutative diagram
$$
\begin{xy}
\xymatrix{
M \ar[r] & C \ar[r]^{c} & M_{\chi_{10}} \\
R \ar[u] \ar[r] & I \ar[u] \ar[r] & R_{\chi_{10}}\ar[u] \\
}\end{xy} \eqno(2)
$$
Note that the ring $M$ of vector-valued modular forms is
not finitely generated, but the ring $C$ of covariants is.
Any covariant defines a meromorphic Siegel modular form the polar locus
of which is either empty or the divisor $\A{1,1}$. 
A covariant, that is, an equivariant section of 
${\rm Sym}^p (V) \otimes \det(V)^q$ with $\lambda=(p+q,q)$, defines a section 
$F$ of 
$$
{\rm Sym}^p({\EE}) \otimes \det{\EE}^{(q-p)/3}= 
{\EE}_{\lambda} \otimes \det{\EE}^{-2d}
$$ 
on ${\mathcal A}_2-{\mathcal A}_{1,1}$.
If $F$ extends to a holomorphic section of some 
${\rm Sym}^j({\EE}) \otimes \det ({\EE})^k$
on ${\mathcal A}_2$, then it extends as a section of that bundle to all of 
$\tilde{\mathcal A}_2$ by the Koecher principle. 
If the section $F$ has an order $-\nu(F)$ along the divisor ${\mathcal A}_{1,1}$
then multiplying $F$ with $\chi_{10}^{\lceil {\nu(F)}/2 \rceil}$ 
makes it into a holomorphic section. Similarly, multiplication with 
$\chi_5^{\nu(F)}$ makes it into a holomorphic modular 
form with character of weight $(p, (q-p)/3+5\nu(F))$, see the next section. 
In particular, we can apply this to the universal binary sextic $f$ given in 
Example \ref{tautologicalf}. 
It thus defines a meromorphic Siegel modular form of weight $(6,-2)$.
\end{section}
%%%%%%%%%%%%%%%%%%%%%%%%%%%%%%%%%%%%%%%%
\begin{section}{The Forms $\chi_{10}$ and $\chi_{6,8}$}\label{theforms}
In this paper two modular forms (and two closely related ones) 
will play a central role. The modular forms are
$\chi_{10}$ and $\chi_{6,8}$ and the related ones $\chi_5$ and $\chi_{6,3}$. 
Up to a normalization the modular form $\chi_{5}$ is defined as the product of
the ten even theta constants and $\chi_{10}=\chi_5^2$ is its square. 
The form $\chi_5$ can be seen as a cusp form of weight $5$ on 
the congruence subgroup $\Gamma_2[2]$ of level $2$ which is invariant
under the alternating group $\mathfrak{A}_6 
\subset \mathfrak{S}_6=\Sp{4}{\ZZ/2\ZZ}$;
alternatively it can be seen as a modular form on $\Gamma_2$ with
character $\epsilon$. (The abelianization of $\Gamma_2$ is 
${\ZZ}/2{\ZZ}$, see \cite[Satz 3]{Klingen}, \cite{Maass}.) The Fourier expansion of $\chi_5$ follows 
from its definition and starts as
$$
\begin{aligned}
\chi_5(\tau)=(1/R-R)\, Q_1Q_2+(-1/R^3-9/R+R^3+9R)(Q_1^3Q_2+Q_1Q_2^3) &\\ 
+(-9/R^5+93/R^{3}-90/R+90R-93R^3+9R^5)\, Q_1^3Q_2^3\, + & \cdots \, ,\\
\end{aligned}
$$
where we write $\tau \in \mathfrak{H}_2$ as 
$$
\tau=(\begin{matrix}
\tau_1 & \tau_{12} \\ \tau_{12} & \tau_2 \\ 
\end{matrix})
$$
and use 
$Q_1=e^{\pi i \tau_1}$, $Q_2=e^{\pi i \tau_2}$ and $R=e^{\pi i \tau_{12}}$. 
The Fourier expansion of $\chi_{10}$ thus starts as follows:
$$
\begin{aligned}
\chi_{10}(\tau)=
(1/r-2+r)q_1q_2-(2/r^2+16/r+36+16r+2r^2)(q_1^2q_2+q_1q_2^2)+ & \\
(1/r^3+36/r^2+99/r-272+99r36r^2+r^3)(q_1^3q_2+q_1q_2^3) + &  \\
(-16/r^3+240/r^2-240/r+32-240r+240r^2-16r^3)q_1^2q_2^2 \, + & \cdots\, , \\
\end{aligned}
$$
where
$q_1=e^{2\pi i \tau_1}$, $q_2=e^{2\pi i \tau_2}$ and $r=e^{2\pi i \tau_{12}}$.

We will be interested in the development of $\chi_5$ and $\chi_{10}$
along the locus $\mathfrak{H}_1\times \mathfrak{H}_1\subset \mathfrak{H}_2$ 
given by $\tau_{12}=0$.
Note that the pullback to ${\mathcal A}_1 \times {\mathcal A}_1$ 
of the normal bundle along ${\mathcal A}_{1,1}$ is the dual of the 
product ${\EE} \boxtimes {\EE}$ of the (pullback of the) Hodge bundles 
of the two factors,
see \cite[p.\ 23]{C-vdG}.

We need the concept of a quasi-modular form and refer to \cite{KZ} and 
to \cite{Zagier}.
For even $k\geq 2$ we will denote the Eisenstein series of weight $k$ 
on $\SL{2}{\ZZ}$ by $e_k$. 
Its Fourier expansion is given by
\[
e_{k}(\tau)=
1-\frac{2k}{B_k}\sum_{n\geqslant 1}\sigma_{k-1}(n) \, q^n.
\]
with $B_k$ the $k$th Bernoulli number and $\sigma_r(n)=\sum_{d|n} d^r$.
For a subgroup $\Gamma$ of 
finite index of the full modular group $\SL{2}{\ZZ}$ we denote by
$M_{*}(\Gamma)=\oplus_k M_k(\Gamma)$ the graded ring of modular forms and
by
$
\widetilde{M_*}(\Gamma)=\oplus_k \widetilde{M_k}(\Gamma)
$
the graded ring of quasi-modular forms on $\Gamma$.
We have a differential operator
$$
D=\frac{1}{2 \pi i}\frac{d}{d\tau}=q \, \frac{d}{dq}
$$
that sends quasi-modular forms to quasi-modular forms.
We refer to \cite{KZ} for the following facts.

\begin{lemma} 
We have
$\widetilde{M_*}(\Gamma)=M_*(\Gamma)\otimes \CC[e_2]$.
Furthermore, for even $k > 0$ we have
$$\widetilde{M_k}(\Gamma)=
\oplus_{0\leqslant i \leqslant k/2} D^i M_{k-2i}(\Gamma)
\oplus \langle D^{k/2-1}e_2
\rangle \, .
$$ 
\end{lemma}

We develop $\chi_{10}$ as a Taylor series in the 
normal direction of $\mathfrak{H}_1 \times \mathfrak{H}_1$ inside $\mathfrak{H}_2$
with coordinate $t=2\pi i \, \tau_{12}$
$$
\chi_{10}(\tau)= \sum_{m=0}^{\infty} \xi_{m} \frac{t^m}{m!} \, ,
\quad \text{\rm with}  \quad 
\xi_m= \frac{\partial^m \, \chi_{10}}{\partial t^m } 
|_{\left(
\begin{smallmatrix}
\tau_1 & 0 \\ 
0 & \tau_2 \\ 
\end{smallmatrix}
\right)}  \, .
$$
Since sections of the Hodge bundle on $\tilde{\mathcal A}_1$ 
are modular forms and since we know that the operator $D$
sends the ring of quasi-modular forms on ${\rm SL}(2,{\ZZ})$ to itself,
it follows that the terms in the Taylor series of $\chi_{10}$ on $\mathfrak{H}_2$
along $\mathfrak{H}_1\times \mathfrak{H}_1$ are quasi-modular forms.
Using the definition in terms of even theta constants we can calculate the expansions
of $\chi_{10}$ and $\chi_5$.

\begin{lemma} We have $\xi_m=0$ for $m$ odd. 
The first three non-zero coefficients of $\chi_{10}=\sum_{m=0}^{\infty} \xi_m \, t^m/m!$ 
are 
$$
\xi_2= 2 \, \Delta \otimes \Delta\, , \qquad
\xi_4= 2 \, \Delta e_2  \otimes \Delta e_2\, ,
$$ 
and
$$
\xi_6=\frac{1}{24} \left(-7 \, \Delta e_4 \otimes
\Delta e_4 +65 \, \Delta e_2^2 \otimes \Delta e_2^2
-5 \, (\Delta e_4 \otimes \Delta e_2 ^2
+\Delta e_2^2\otimes \Delta e_4) \right) \, .
$$
Here we use the shorthand  $f \otimes g$ instead of $f(\tau_1)\otimes g(\tau_2)$.
\end{lemma}

Similarly, for the development of $\chi_5$ as a Taylor series in 
the normal direction we need the square root of $\Delta$  
which is the modular form
$$
\delta= q^{1/2} \prod_{n=1}^{\infty} (1-q^n)^{12}
$$
of weight $6$ on the congruence subgroup $\Gamma_1(2)$.
If we use $s=\pi i \tau_{12}$ as the normal coordinate 
then we find the development
$$
\chi_5(\tau)= 
-2\, \delta \otimes \delta \, s 
- \frac{1}{3}  \delta e_2 \otimes \delta e_2 \, s^3 +\cdots \, .
$$

The second form that plays a central role is the form $\chi_{6,8}$
(and its relative $\chi_{6,3}$).
We defined this vector-valued Siegel modular form of weight 
$(6,3)$ with character $\epsilon$ on $\Gamma_2$ in the paper \cite[Example 16.1]{C-vdG-G}
using the gradients
of the six odd theta series $\vartheta_{m_i}(\tau,z)$ ($i=1,\ldots,6)$ 
with characteristics in degree $2$:
$$
\chi_{6,3}=c\cdot \Sym^6(G_1,G_2,G_3,G_4,G_5,G_6)
\quad
\text{with}
\quad
G_i(\tau)=
\left(
\begin{smallmatrix}
{\partial \vartheta_{m_i}}/{\partial z_1}\\
{\partial \vartheta_{m_i}}/{\partial z_2}
\end{smallmatrix}
\right)(\tau,(0,0)) \, ; \eqno(3)
$$
here the constant $c\in \CC^*$ is chosen such that the Fourier expansion of
$\chi_{6,3}$ starts as follows
$$
\chi_{6,3}(\tau)=
\left(
\begin{smallmatrix}
0\\
0\\
(R-R^{-1})\\
(2R+2R^{-1})\\
(R-R^{-1})\\
0\\
0
\end{smallmatrix}
\right)Q_1Q_2+\cdots \, , \eqno(4)
$$
where $Q_1=e^{i\pi \tau_1}$, $Q_2=e^{i\pi \tau_2}$ and $R=e^{i\pi \tau_{12}}$.
The form $\chi_{6,3}$ is a modular form of weight $(6,3)$ with character
$\epsilon$ on $\Gamma_2$. 
Alternatively, it can be viewed as a form on the level $2$ 
principal congruence subgroup $\Gamma_2[2]$ that is invariant under  
the action of the alternating group 
$\mathfrak{A}_6$. 
We define $\chi_{6,8}\in S_{6,8}(\Gamma_2)$ as the product 
$\chi_{6,3} \chi_5$.

Using the definition with the gradients 
we find the Taylor expansion of $\chi_{6,3}$
in the normal direction along $\mathfrak{H}_1 \times \mathfrak{H}_1$ with $s=\pi i \tau_{12}$
as coordinate:
$$
\chi_{6,3}(\tau)=
\left(
\begin{smallmatrix}
0\\
0\\
0\\
4\, \delta \otimes \delta \\
0\\
0\\
0
\end{smallmatrix}
\right)+
\left(
\begin{smallmatrix}
0\\
0\\
2\, e_2\delta \otimes \delta  \\
0\\
2\, \delta\otimes e_2\delta  \\
0\\
0
\end{smallmatrix}
\right)
\, s +
\left(
\begin{smallmatrix}
0\\
\frac{2}{3}\, (e_2^2-e_4)\delta \otimes \delta \\
0 \\
4\, e_2\delta \otimes e_2\delta \\
0 \\
\frac{2}{3}\, \delta \otimes (e_2^2-e_4)\delta \\
0
\end{smallmatrix}
\right) \, 
\frac{s^2}{2 !}+\cdots \, .
$$

This gives then the Taylor development of $\chi_{6,8}$ in the normal direction
with coordinate $t=2 \pi i \tau_{12}$:
$$
\chi_{6,8}(\tau)=
-
\left(
\begin{smallmatrix}
0\\
0\\
0\\
4\, \Delta \otimes \Delta \\
0\\
0\\
0
\end{smallmatrix}
\right) \, t -
\left(
\begin{smallmatrix}
0\\
0\\
e_2\Delta \otimes \Delta \\
0\\
\Delta \otimes e_2\Delta \\
0\\
0
\end{smallmatrix}
\right)\, t^2-
\left(
\begin{smallmatrix}
0\\
2(e_2^2-e_4)\Delta \otimes \Delta \\
0\\
16 e_2\Delta \otimes e_2\Delta \\
0\\
2\Delta \otimes (e_2^2-e_4)\Delta \\
0
\end{smallmatrix}
\right) \, \frac{t^3}{3} +\cdots \, .
$$
The modular form $\chi_{6,8}$ was first given by 
Ibukiyama using theta series with pluri-harmonic coefficients, 
see \cite{Ibukiyama}. It is the first vector-valued (non scalar-valued)
modular form if one orders them according Deligne weight ($=j+2k-3$).
\end{section}
%%%%%%%%%%%%%%%%%%%%%%%%%%%%%%%%%%%%%%%%%%
\begin{section}{Modular Forms Associated to Covariants}\label{MFassociatedtoCovariants}
Recall from Section \ref{covariantsandSMF} the inclusions
$$
M \longrightarrow C {\buildrel c \over \longrightarrow} M_{\chi_{10}} \, . 
$$
We wish to describe the map $c$ explicitly.
For this we consider the meromorphic modular form $\chi_{6,3}/\chi_{5}$ 
of weight $(6,-2)$ (without character). 
It is holomorphic on $\tA{2}-\A{1,1}$, hence defines by Lemma \ref{pullbackE} 
under pullback a covariant, a section of $\Sym^6(V)$.
Up to a non-zero multiplicative scalar $r$ we thus find the tautological~$f$
(see Example \ref{tautologicalf}), hence
$$
r \, c(f)=\frac{\chi_{6,3}}{\chi_5} \, .
$$
We will discard the factor $r$ and assume that $r=1$ by normalizing $f$ appropriately.
But every covariant can be constructed from $f$. In fact, let $h$ be a
covariant associated to the irreducible representation 
$A[p+q,q]$ in $\Sym^d(\Sym^6(V))$. It 
can be viewed as a form of degree $d$ in $a_0,\ldots,a_6$, the coefficients of $f$, 
and of degree $p$ (resp.\ $q$) 
in the coordinates $x_1$ and $x_2$, or simply as a
vector of length $p+1$ with entries which are polynomials of degree $d$ in 
$a_0,\ldots,a_6$. 

Now write the Fourier expansion of $\chi_{6,3}$ as a vector 
$(\alpha_0,6\, \alpha_2,15\, \alpha_2, 20\, \alpha_3, \ldots, \alpha_6)^t$, 
or more symbolically as
$$
\chi_{6,3}= \sum_{i=0}^6  \binom{6}{i} \alpha_i \, X^{6-i}Y^i\, ,
$$
where $\alpha_i$ is a Taylor series in $Q_1,Q_2$ and $R$ and $R^{-1}$,
and the $X^{6-i}Y^i$ indicate the coordinate places. 
Define a map
$$
\gamma: C \to M,\qquad  h \mapsto F_h,
$$
where $F_h$ is obtained by substituting $\alpha_i$ for $a_i$ in $h$.
When $h$ is viewed as a vector of length $p+1$  this 
gives the $p+1$ entries of a holomorphic vector-valued modular form $F_h$ 
of weight $(p,q+3d)$. In particular, we have
$$
c(h)= \gamma(h)/\chi_5^d
$$ 
and we see that the order of $c(h)$ along $\mathfrak{H}_1^2$ is $\geq -5d$.
Using the expansion of $\chi_{6,3}$ and $\chi_5$ along $\mathfrak{H}_1^2$ 
given in Section \ref{theforms} 
we can calculate the order $\nu(F_h)$ of $F_h$ along $\mathfrak{H}_1^2$. 
Division by $\chi_5^{\nu(F_h)}$ yields a holomorphic modular form
on $\Gamma_2$ (with or without a character) of 
weight $(p,q+3d-5\nu(F_h))$.

\begin{conclusion}
Every vector-valued modular form of given weight 
on $\Gamma_2$ can be constructed from a covariant
by substituting the Fourier coefficients of $\chi_{6,3}$ and by multiplying with
an appropriate power of $\chi_5$. The same holds for modular forms on $\Gamma_2$
with a character.
\end{conclusion}

We shall see in the next sections that this gives a very effective way of
constructing the Fourier expansions of vector-valued modular forms of degree $2$.
In fact, we can calculate the Fourier expansion of $\chi_{6,3}$ and $\chi_5$ 
easily as these are given in terms
of theta functions.
%\smallskip

\begin{remark}
The central role that $\chi_{6,3}$ plays can be explained as follows.
It is well-known that for a smooth curve $C$ of genus $2$ there are six symmetric 
translates of the theta divisor, each isomorphic to $C$,  
in the Jacobian ${\rm Jac}(C)$ passing through
the origin. The tangents to these six curves give six points in the projectivized
tangent space at the origin, that is, six points on ${\PP}^1$. That is the way 
to retrieve the sextic as was remarked by Bolza, see \cite[p.\ 481]{Bolza}. 
In terms of odd theta characteristics this means that
if we write the Taylor expansion of the odd theta functions as
$$
\vartheta_{m_i}=c_{i,1}\, z_1+c_{i,2}\, z_2 + \text{\rm higher order terms,}
$$
then we get the form $\chi_{6,3}$ up to a normalization 
as the product of these linear terms
$$
\chi_{6,3} \sim \prod_{i=1}^6 (c_{i,1} z_1 + c_{i,2} z_2)= c_{1,1}\cdots c_{6,1} z_1^6
+\ldots + c_{1,2}\cdots c_{6,2}z_2^6 \, .
$$
\end{remark}
  
\end{section}

%%%%%%%%%%%%%%%%%%%%%%%%%%%%%%%%%%%%%%%%
\begin{section}{Construction of Modular Forms}
In this section we give examples of constructions of vector-valued and
scalar-valued modular forms using covariants. Our point of departure
is the Fourier series of $\chi_{6,3}$ that we calculated as in (4)
as a series in $Q_1$ and $Q_2$
up to $Q_1^aQ_2^b$ with $a+b \leq 170$.

\medskip

\centerline{\fbox{\parbox{1.0cm}{\bf $d=2$}}}

\smallskip

We start with $d=2$, see Example \ref{example-d2}.
Covariants associated to the isotypical decomposition
$$
\Sym^2(\Sym^6(V))=A[12,0]+A[10,2]+A[8,4]+A[6,6]
$$
provide by the construction given above modular forms in
$S_{12,6}$, $S_{8,8}$, $S_{4,10}$ and $S_{0,12}$.
Notice that all the latter spaces are one dimensional.
The modular form in $S_{12,6}$ is the square of $\chi_{6,3}$
and this gives us its Fourier expansion immediately.

The covariant corresponding to $s[10,2]$ is the Hessian
and we thus find a modular form $\chi_{8,8} \in S_{8,8}$ with coordinates
(symmetric in the sense that interchanging 
$\alpha_i$ and  $\alpha_{6-i}$ reverses the vector)
$$
\begin{aligned}
(\alpha_0 \alpha_2-\alpha_1^2, 4\alpha_0\alpha_3-4\alpha_1\alpha_2,
6\alpha_0\alpha_4+4\alpha_1\alpha_3-10 \alpha_2^2,
4\alpha_0\alpha_5+16\alpha_1\alpha4-20\alpha_2\alpha_3, &\\
\alpha_0\alpha_6+14\alpha_1\alpha_5+5\alpha_2\alpha_4-20\alpha_3^2,
 \ldots,
\alpha_4\alpha_6 -\alpha_5^2) & \\
\end{aligned}
$$
which we normalize such that its Fourier expansion starts with:
$$
\chi_{8,8}(\tau)=
\left(
\begin{smallmatrix}
0\\
0\\
(r-2+r^{-1})\\
3(r-r^{-1})\\
(4\, r+10+4\, r^{-1})\\
3(r-r^{-1})\\
(r-2+r^{-1})\\
0\\
0
\end{smallmatrix}
\right)q_1q_2+\ldots \, ,
$$
where as before
$q_1=e^{2 i\pi \tau_1}$, 
$q_2=e^{2 i\pi \tau_2}$ and $r=e^{2 i\pi \tau_{12}}$.
By restriction along $\mathfrak{H}_1 \times \mathfrak{H}_1$ 
we find a non-zero multiple of the transpose of
$
(0,\ldots,0,\Delta \otimes \Delta,0,\ldots,0)
$
which shows that $\chi_{8,8}$ is not divisible by $\chi_5$.
Similarly, the covariant yielding a modular form in $S_{4,10}$ has coordinates
that are reversed  under interchanging $\alpha_i$ and $\alpha_{6-i}$
$$
(\alpha_0\alpha_4-4 \alpha_1\alpha_3+3\alpha_2^2, 2\alpha_0\alpha_5-
6\alpha_1\alpha_4+4\alpha_2\alpha_3, \alpha_6\alpha_9-9\alpha_2\alpha_4
+8\alpha_3^2, \ldots, \alpha_2\alpha_6-4\alpha_3\alpha_5+3\alpha_4^2)
$$ 
We normalize the resulting form $\chi_{4,10}$ such that 
$$
\chi_{4,10}(\tau)=
\left(
\begin{smallmatrix}
(r-2+r^{-1})\\
2(r-r^{-1})\\
3(r+6+r^{-1})\\
2(r-r^{-1})\\
(r-2+r^{-1})
\end{smallmatrix}
\right)q_1q_2+\cdots \, .
$$
Its restriction to $\mathfrak{H}_1\times \mathfrak{H}_1$ is 
a non-zero multiple of 
$(0,0,\Delta \times \Delta,0,0)^t$, and we thus cannot divide by $\chi_5$. 
Finally, we get a non-zero form in $S_{0,12}$ by taking
$$
\chi_{12}=  \alpha_0\alpha_6 -6 \, \alpha_1\alpha_5 +15 \, 
\alpha_2 \alpha_4 - 10\,  \alpha_3^2 
$$
which we normalize so that it starts by $\chi_{12}=(r+10+1/r)q_1q_2+ \cdots$.
One immediately gets the Fourier expansion of $\chi_{12}$.

In all these cases we checked the Fourier expansion by calculating 
Hecke eigenvalues for the Hecke operators $T(p)$ for primes $p\leq 23$
and checked that these fit with the values provided by \cite{FvdG} and 
\cite{BFvdG}. We give the Hecke eigenvalues for $\chi_{8,8}$ and $\chi_{12,6}$ in a table in Section \ref{tables}.

\smallskip

\centerline{\fbox{\parbox{1.0cm}{\bf $d=3$}}}

\smallskip
In this case we have the isotypical decomposition
$$
\Sym^3(\Sym^6(V))= A[18, 0] + A[16, 2] + A[15, 3] + A[14, 4] + A[13, 5] + 2\, A[12, 6] + A[10, 8]
$$
which by the procedure of the last section leads to modular forms of weights
$(18,9)$, $(14,11)$, $(12,12)$, $(10,13)$, $(8, 14)$, $(6,15)$ and $(2,17)$
on the group $\Gamma_2$ with a character.
If these forms are divisible by $\chi_5$ or $\chi_{10}$
we will find forms of smaller weight.
Note that in the case of weight $(6,15)$ we find
a $2$-dimensional space of modular forms.
In the cases of weight $(18,9)$, $(14,11)$, $(10,13)$, $(2,17)$ 
the restriction of the corresponding form to
$\mathfrak{H}_1 \times \mathfrak{H}_1$ is a non-zero multiple of the transpose of
a vector of the form
$$
(0,\ldots,0,\Delta \delta \otimes \Delta \delta,0,\ldots,0)
$$
and thus gives non-zero modular forms and these forms are not divisible by $\chi_5$.
We analyze the remaining three cases. In the case of weight $(12,12)$ our form vanishes
with order $2$ on $\mathfrak{H}_1^2$; dividing by $\chi_{10}$ 
yields a form $\chi_{12,2}$ on $\Gamma_2$ with character whose Fourier expansion starts with
$$
\chi_{12,2}(\tau)=
\left(
\begin{smallmatrix}
0\\
0\\
0\\
2(R-R^{-1})\\
9(R+R^{-1})\\
12(R-R^{-1})\\
0\\
-12(R-R^{-1})\\
-9(R+R^{-1})\\
-2(R-R^{-1})\\
0\\
0\\
0
\end{smallmatrix}
\right)Q_1Q_2+\cdots \, .
$$
Note that its square ${\rm Sym}^2 \chi_{12,2}$ is the generator of the space $S_{24,4}$.
Similarly, in the case of weight $(8,14)$ 
our form vanishes with order $2$ along $\mathfrak{H}_1^2$
and by dividing by $\chi_{10}$ we find 
a form of weight $(8,4)$ with character.

The representation $A[12,6]$ occurs with multiplicity $2$ in $\Sym^3(\Sym^6(V))$
and we thus can find two linearly independent covariants, say $h_1$ and $h_2$.
The general linear combination of these two gives a modular form that 
upon restriction to
$\mathfrak{H}_1^2$ yields a non-zero multiple of the transpose of
$$
[0,0,0,\eta^{12}\, \Delta(\tau_1)\otimes \eta^{12}\,\Delta(\tau_2),0,0,0] \, .
$$
The linear combination that vanishes along $\mathfrak{H}_1^2$ vanishes with
multiplicity $2$ and by division by $\chi_{10}$ we find a cusp form with
character of weight $(6,5)$ the Fourier expansion of which starts with
$$
\chi_{6,5}(\tau)=
\left(
\begin{smallmatrix}
2(R-R^{-1})\\
6(R+R^{-1})\\
5(R-R^{-1})\\
0\\
5(R-R^{-1})\\
6(R+R^{-1})\\
2(R-R^{-1})
\end{smallmatrix}
\right)Q_1Q_2+\cdots \, .
$$
\end{section}
%%%%%%%%%%%%%%%%%%%%%%%%%%%%%%%%%%%%%%%
\begin{section}{Further Examples}
Here we look at a few cases with $d=4$ and $d=5$.

\smallskip

\centerline{\fbox{\parbox{1.0cm}{\bf $d=4$}}}

\smallskip

We have the isotypical decomposition
$$
\begin{aligned}
\Sym^4(\Sym^6(V))= A[24, 0] + A[22, 2] + A[21, 3] + 2\, A[20, 4]
+ A[19, 5] + 3\, A[18, 6] + &\\
A[17,7]+3\, A[16,8] +A[15,9]+2\, A[14,10]+2\, A[12, 12]\, . & \\
\end{aligned}
$$
We then use the covariants to construct modular forms from $\chi_{6,3}$.
By calculating the behavior along $\mathfrak{H}_1^2$
we can deduce the orders of vanishing along $\mathfrak{H}_1^2$.

We claim that the following table gives
the orders of the corresponding modular forms along $\mathfrak{H}_1^2$.
If the multiplicity
$\mu_{m,n}$ is greater than $1$ we list the orders of vanishing occurring in
the $\mu_{m,n}$-dimensional space of covariants.

\begin{footnotesize}\smallskip
\vbox{
\bigskip\centerline{\def\quad{\hskip 0.6em\relax}
\def\quod{\hskip 0.5em\relax }
\vbox{\offinterlineskip
\hrule
\halign{&\vrule#&\strut\quod\hfil#\quad\cr
height2pt&\omit&&\omit&&\omit&&\omit&&\omit&\cr
%\noalign{\hrule}
& $[m,n]$  && $\mu_{m,n}$ && weight && order &\cr
\noalign{\hrule}
& $[24,0]$ && $1$ && $(24,12)$ && $0$ & \cr
& $[22,2]$ && $1$ && $(20,14)$ && $0$ &\cr
& $[21,3]$ && $1$ && $(18,15)$ && $2$ &\cr
& $[20,4]$ && $2$ && $(16,16)$ && $0,2$ &\cr
& $[19,5]$ && $1$ && $(14,17)$ && $2$ &\cr
& $[18,6]$ && $3$ && $(12,18)$ && $0,2,3$ &\cr
& $[17,7]$ && $1$ && $(10,19)$ && $2$ &\cr
& $[16,8]$ && $3$ && $(8,20)$ && $0,2,3$ &\cr
& $[15,9]$ && $1$ && $(6,21)$ && $2$ &\cr
& $[14,10]$ && $2$ && $(4,22)$ && $0,2$ &\cr
& $[12,12]$ && $2$ && $(0,24)$ && $2,4$ &\cr
} \hrule}
}}
\end{footnotesize}

From this table one can read off what modular forms 
can be constructed. For example,
for the representation $A[18,6]$ (resp.\ $A[16,8]$) we find 
multiplicity $3$, hence a $3$-dimensional space of cusp forms
of weight $(12,18)$ (resp.\ weight $(8,20)$). 
We can calculate generating modular forms using the expressions for
the covariants and observe
that in both cases there is a $2$-dimensional subspace of forms vanishing with
order $2$ along $\mathfrak{H}_1^2$. Dividing these forms by $\chi_{10}$ yields a
$2$-dimensional space of cusp forms of weight $(12,8)$ (resp.\ $(8,10)$). 
By analyzing their behavior
along $\mathfrak{H}_1^2$ we see that there is a $1$-dimensional
subspace vanishing along $\mathfrak{H}_1^2$ and division by
$\chi_5$ leads to a cusp form of weight $(12,3)$ (resp.\ $(8,5)$) with character.
Alternatively, from the vanishing/non-vanishing of spaces of cusp forms
(with/without a character) one can read off the orders of vanishing.

We give the details for $A[16,8]$. We find a $3$-dimensional space
of cusp forms $\gamma(h)$ in $M$ as $h$ ranges over the space of covariants 
associated to $d=4$ and $A[16,8]$. The generic element yields
$[0,\ldots,0,\Delta^2\otimes\Delta^2,0,\ldots,0]$ when restricted to
$\mathfrak{H}_1^2$ and there is a $2$-dimensional space of forms of weight
$(8,10)$ vanishing at least doubly on $\mathfrak{H}_1^2$. Division by
$\chi_{10}$ gives a $2$-dimensional space of forms of weight $(8,10)$ 
generated by $G_1$ and $G_2$ with Fourier expansions
$$
G_1(\tau)=\left(
\begin{smallmatrix}
(r-2+r^{-1})\\
4(r-r^{-1})\\
(9\,r+34+9\,r^{-1})\\
13(r-r^{-1})\\
15(r-2+r^{-1})\\
13(r-r^{-1})\\
(9\,r+34+9\,r^{-1})\\
4(r-r^{-1})\\
(r-2+r^{-1})
\end{smallmatrix}
\right)q_1q_2+\cdots \, ,
\quad
G_2(\tau)=\left(
\begin{smallmatrix}
3(r-2+r^{-1})\\
12(r-r^{-1})\\
2(11\,r+26+11\,r^{-1})\\
24(r-r^{-1})\\
25(r-2+r^{-1})\\
24(r-r^{-1})\\
2(11\,r+26+11\,r^{-1})\\
12(r-r^{-1})\\
3(r-2+r^{-1})
\end{smallmatrix}
\right)q_1q_2+\cdots \, .
$$
The action of the Hecke operator $T_2$ on $G_1$ and $G_2$ is as follows:
$$
G_1\vert_{8,10} T_2=8160\, G_1-4080\, G_2,\quad
G_2\vert_{8,10} T_2=-2880\, G_1-720\, G_2.
$$
We get two Hecke eigenforms in $S_{8,10}$, namely 
$\chi_{8,10}^{(\pm)}=(37\pm\sqrt{2185})G_1-34\,G_2$, 
on which the Hecke operator $T_2$ acts with
eigenvalues $(3720 \pm 120\sqrt{2185})$.

For the case of $A[18,6]$, leading to modular forms of weight $(12,18)$ and
after division by $\chi_{10}$ to a $2$-dimensional space of cusp forms of
weight $(12,8)$, the story is very similar  and we can caluculate the Fourier
expansions and Hecke eigenvalues as well. We list these eigenvalues for weight
$(8,10)$ and $(12,8)$ in the following table.

\begin{footnotesize}
\smallskip
\vbox{
\bigskip\centerline{\def\quad{\hskip 0.6em\relax}
\def\quod{\hskip 0.5em\relax }
\vbox{\offinterlineskip
\hrule
\halign{&\vrule#&\strut\quod\hfil#\quad\cr
height2pt&\omit&&\omit &&\omit &\cr
%\noalign{\hrule}
& $q$ &&
$\lambda_q(\chi_{8,10}^{\pm})$  && $\lambda_{q}(\chi_{12,8}^{\pm})$ &\cr
\noalign{\hrule}
height2pt&\omit&&\omit &&\omit &\cr
%& && && &\cr
& $2$ && $3720 \pm 120 \sqrt{2185}$ && $-768 \pm 192\sqrt{1381}$  &\cr
& $3$ && $674280 \mp 18720\sqrt{2185}$  && $86616 \pm 20736\sqrt {1381}$ &\cr
& $4$ && $-945536 \mp 28800\sqrt{2185}$ && $790528 \mp 147456\sqrt {1381}$  &\cr
& 5 && $-70706100 \pm 8188800\sqrt{2185}$  && $362491500 \mp 5145600\sqrt {1381}$ &\cr
& 7 && $-11441284400 \mp 446644800\sqrt{2185}$  && $14252364592 \mp 459468288\sqrt {1381}$ &\cr
} \hrule}
}}
\end{footnotesize}

The restriction of $G_1$ (resp.\ $G_2$) to $\mathfrak{H}_1^2$ is $52$ (resp.\
$96$) times the transpose of 
$$
[0,0,e_{4}\Delta \otimes\Delta,0,0,0,\Delta\otimes e_{4}\Delta,0,0]
$$
so $G_1/52-G_2/96$ 
vanishes  on $\mathfrak{H}_1^2$ and dividing it by $\chi_{5}$
gives the unique cusp form of weight $(8,5)$ on $\Gamma_2[2]$ which
is $\mathfrak{S}_6$-anti-invariant.

Finally, in the case of $A[12,12]$ we find a $2$-dimensional space of cusp forms
of weight $(0,24)$. Every modular form in this space restricts to a multiple
of $\Delta^2\otimes \Delta^2$ and there is a linear combination that vanishes
on $\mathfrak{H}_1^2$ and it is a multiple of $E_4 \, \chi_{10}^2$. Thus we see
that by dividing by $\chi_{10}^2$ we get the Eisenstein series of weight $4$. 
\bigskip

\smallskip

\centerline{\fbox{\parbox{1.0cm}{\bf $d=5$}}}

\smallskip

We list the irreducible representations occurring in the isotypical decomposition
for $d=5$ in a table together with their multiplicities, the weight of the
corresponding modular forms and the orders of vanishing along $\mathfrak{H}_1^2$.

\begin{footnotesize}\smallskip
\vbox{
\bigskip\centerline{\def\quad{\hskip 0.6em\relax}
\def\quod{\hskip 0.5em\relax }
\vbox{\offinterlineskip
\hrule
\halign{&\vrule#&\strut\quod\hfil#\quad\cr
height2pt&\omit&&\omit&&\omit&&\omit&&\omit&&\omit&&\omit&&\omit&&\omit&\cr
%\noalign{\hrule}
& $[m,n]$  && $\mu_{m,n}$ && weight && order && && $[m,n]$ && $\mu_{m,n}$ && weight && order &\cr
\noalign{\hrule}
& $[30,0]$ && $1$ && $(30,15)$ && $0$ && && $[22,8]$ && $4$ && $(14,23)$ && $0,2,3,4$ & \cr
& $[28,2]$ && $1$ && $(26,17)$ && $0$ && && $[21,9]$ && $3$ && $(12,24)$ && $2,3,4$ & \cr
& $[27,3]$ && $1$ && $(24,18)$ && $2$ && && $[20,10]$ && $4$ && $(10,25)$ && $0,2,3,4$ & \cr
& $[26,4]$ && $2$ && $(22,19)$ && $0,2$ && && $[19,11]$ && $2$ && $(8,26)$ && $2,3$ & \cr
& $[25,5]$ && $2$ && $(20,20)$ && $2,3$ && && $[18,12]$ && $4$ && $(6,27)$ && $0,2,3,4$ & \cr
& $[24,6]$ && $3$ && $(18,21)$ && $0,2,3$ && && $[17,13]$ && $1$ && $(4,28)$ && $2$ & \cr 
& $[23,7]$ && $2$ && $(16,22)$ && $2,3$ && && $[16,14]$ && $2$ && $(2,29)$ && $0,4$ & \cr
} \hrule}
}}
\end{footnotesize}

We give the details for the case of $\lambda=(18,12)$, where the $4$-dimensional
space of covariants produces modular forms of weight $(6,27)$ with
character. We find a 
$2$-dimensional subspace of modular forms vanishing with multiplicity $\geq 3$,
and dividing these by $\chi_5^3$ produces a basis of $S_{6,12}$ which is of 
dimension $2$. The two basis elements $G_1$ and $G_2$ have Fourier expansions starting with
$$
G_1(\tau)=\left(
\begin{smallmatrix}
2(r+10+r^{-1})\\
6(r-r^{-1})\\
33(r-2+r^{-1})\\
56(r-r^{-1})\\
33(r-2+r^{-1})\\
6(r-r^{-1})\\
2(r+10+r^{-1})
\end{smallmatrix}
\right)q_1q_2+\cdots \, ,
\quad
G_2(\tau)=\left(
\begin{smallmatrix}
0\\
0\\
(r-2+r^{-1})\\
2(r-r^{-1})\\
(r-2+r^{-1})\\
0\\
0
\end{smallmatrix}
\right)q_1q_2+\cdots \, .
$$
One calculates the action of the Hecke operator $T_2$ and finds
$$
G_1\vert_{6,12} T_2=-2592\, G_1+254016\, G_2,\quad
G_2\vert_{6,12} T_2=-480\, G_1+24960\, G_2 
$$
and obtains eigenforms 
$\chi_{6,12}^{(\pm)}=(-41\pm \sqrt{601})\, G_1+756\,G_2$.

In a similar way one finds a basis for the following $2$-dimensional
spaces of cusp forms $S_{14,8}$, $S_{12,9}$ and $S_{10,10}$.
One thus can calculate Hecke eigenvalues from the Fourier expansions.
We give the results in the tables in Section \ref{tables}.

The cusp form $G_2$ vanishes identically along
$\mathfrak{H}_1^2$ and $G_2/\chi_{5}$ generates
the $1$-dimensional space $S_{6,7}(\Gamma_2,\epsilon)$.

\end{section}

%%%%%%%%%%%%%%%%%%%%%%%%%%%%%%%%%%%%%%%%
\begin{section}{A Final Example: The Space $S_{4,16}$}
In this section we illustrate the effectivity of our approach and
show how one can use covariants to construct a basis
for the $3$-dimensional space $S_{4,16}$ and use this to calculate 
eigenvalues for the Hecke operators.

In the decomposition of $\Sym^8(\Sym^6(V))$ the irreducible representation
$A[26,22]$ of $\GL{2}{\CC}$ occurs with multiplicity $7$. By the construction of
Section \ref{MFassociatedtoCovariants} this leads to a $7$-dimensional subspace of $S_{4,46}$
vanishing with multiplicity $\geq 7$ at the divisor at infinity.
The restriction to $\mathfrak{H}_1^2$ 
of the cusp forms associated to the corresponding covariants are  all
multiples of the transpose of $[0,0,\Delta^4\otimes \Delta^4,0,0]$.
In fact, we find a $6$-dimensional subspace of forms vanishing with order $\geq 2$
along $\mathfrak{H}_1^2$. Dividing these forms by $\chi_{10}$ leads to 
a $6$-dimensional space of cusp forms
of weight $(4,36)$. These forms restrict to multiples of 
the transpose of $[e_4\, \Delta^3 \otimes \Delta^3, 0, 0, 0, \Delta^3 \otimes 
e_4 \, \Delta^3]$.  Thus there is a $5$-dimensional subspace of forms of weight
$(4,36)$ vanishing along $\mathfrak{H}_1^2$ and these all vanish 
with multiplicity $\geq 2$ along $\mathfrak{H}_1^2$.
We divide again by $\chi_{10}$ to get a $5$-dimensional space of cusp forms
of weight $(4,26)$. All the forms in this space 
restrict to multiples of the transpose of
$[0,0,e_4\Delta^2 \otimes e_4 \Delta^2,0,0]$ along $\mathfrak{H}^2_1$.
This leads to a $4$-dimensional space of cusp forms vanishing with multiplicity
$\geq 1$ along $\mathfrak{H}_1^2$ and we divide these by $\chi_5$ to find
a $4$-dimensional space of cusp forms of weight $(4,21)$ with character. 
Since all these forms restrict to a multiple of the transpose of
$[0, e_6 \delta \Delta \otimes e_4\delta \Delta, 0, e_4\delta\Delta \otimes e_6 \delta\Delta,0]$
it contains a $3$-dimensional
subspace of cusp forms vanishing along $\mathfrak{H}_1^2$. Division by
$\chi_5$ leads to a $3$-dimensional space of cusp forms of weight $(4,16)$.
We normalize these
forms such that their Fourier expansion starts with:
\[
E_1(\tau)=\left(
\begin{smallmatrix}
(r+10+r^{-1})\\
2(r-r^{-1})\\
3(r-2+r^{-1})\\
2(r-r^{-1})\\
(r+10+r^{-1})
\end{smallmatrix}
\right)q_1q_2+
\left(
\begin{smallmatrix}
2(7\, r^2+308\, r+2106+308\, r^{-1}+7\, r^{-2})\\
4(7\, r^2+214\, r-214\, r^{-1}-7\, r^{-2})\\
2(27\, r^2+504\, r-1062+504\, r^{-1}+27\, r^{-2})\\
4(10\, r^2+76\, r-76\, r^{-1}-10\, r^{-2})\\
2(5\, r^2+76\, r+1134+76\, r^{-1}+5\, r^{-2})
\end{smallmatrix}
\right)q_1^2q_2+
\ldots
\]
\[
E_2(\tau)=\left(
\begin{smallmatrix}
(r+10+r^{-1})\\
2(r-r^{-1})\\
3(r-14+r^{-1})\\
2(r-r^{-1})\\
(r+10+r^{-1})
\end{smallmatrix}
\right)q_1q_2+
\left(
\begin{smallmatrix}
2(-11\, r^2+596\, r+1566+596\, r^{-1}-11\, r^{-2})\\
4(11\, r^2+1006\, r-1006\, r^{-1}-11\, r^{-2})\\
2(9\, r^2+2088\, r+5310+2008\, r^{-1}+9\, r^{-2})\\
4(10\, r^2+76\, r-76\, r^{-1}-10\, r^{-2})\\
2(5\, r^2+76\, r+1134+76\, r^{-1}+5\, r^{-2})
\end{smallmatrix}
\right)q_1^2q_2+
\ldots
\]
\[
E_3(\tau)=\left(
\begin{smallmatrix}
(5\,r+104+5\, r^{-1})\\
10(r-r^{-1})\\
(15\, r-138+15\, r^{-1})\\
10(r-r^{-1})\\
(5\,r+104+5\, r^{-1})
\end{smallmatrix}
\right)q_1q_2+
\left(
\begin{smallmatrix}
2(35\, r^2+3808\, r+18306+3808\, r^{-1}+35\, r^{-2})\\
4(35\, r^2+3338\, r-3338\, r^{-1}-35\, r^{-2})\\
2(243\, r^2+8568\, r+10890+8568\, r^{-1}+243\, r^{-2})\\
4(104\, r^2+1892\, r-1892\, r^{-1}-104\, r^{-2})\\
8(13\, r^2+473\, r+2106+473\, r^{-1}+13\, r^{-2})
\end{smallmatrix}
\right)q_1^2q_2+
\ldots
\]
The restrictions of the $e_i$ along the diagonal are the transposes of
$12 \, [e_4^2\Delta\otimes e_4\Delta,0,0,0,e_4\Delta \otimes e_4^2\Delta]$,
$12 \, [e_4^2\Delta\otimes e_4\Delta,0,-3e_6\Delta\otimes e_6\Delta,0, e_4^2\Delta\otimes e_4\Delta]$ and $6\, [19 e_4^2\Delta\otimes e_4\Delta,0,-18 e_6\Delta\otimes e_6\Delta,0, 19 e_4^2\Delta\otimes e_4\Delta]$. So the orders of vanishing along $\mathfrak{H}_1^2$ 
in the $7$-dimensional space  
of cusp forms of weight $(4,46)$ that is the image under the map $\gamma$
of the covariants given by $A[26,22]$ are $\{ 0, 2, 4, 5, 6, 7 \} $.
Finally we remark that the 
form $(E_1/12+E_2/26-E_3/78)/\chi_{5}$ gives a 
cusp form that generates the space
of cusp forms of weight $(6,11)$ with character.

Using the Fourier expansions that we got we can calculate Hecke eigenvalues.
Let $\alpha$ be a root of the irreducible polynomial
$x^3-1042\, x^2+215915\, x +6800500$ in ${\QQ}[x]$.
This field has discriminant $43803704$ and
a  basis of the ring of integers of $K={\QQ}(\alpha)$  
is given by $1$, $\alpha$ and 
$\beta=(1/16785)\alpha^2+(1062/1865)\, \alpha +3007/3357$.
Then a Hecke eigenform $\chi_{4,16}$  with integral Fourier coefficients
is given by
$$
\begin{aligned}
\chi_{4,16}=
(405064\,\alpha-662955\,\beta+1366955)\, E_1+435\, \alpha\, E_2+ &\\
(-61595\,\alpha+100650\,\beta-210400)\,E_3 \, . &\\
\end{aligned}
$$
The Hecke eigenvalues generate the totally real field ${\QQ}(\alpha)$ 
of degree $3$ over ${\QQ}$ of discriminant
$43803704$. We give some Hecke eigenvalues 
for this form.

\begin{footnotesize}
\smallskip
\vbox{
\bigskip\centerline{\def\quad{\hskip 0.6em\relax}
\def\quod{\hskip 0.5em\relax }
\vbox{\offinterlineskip
\hrule
\halign{&\vrule#&\strut\quod\hfil#\quad\cr
height2pt&\omit&&\omit &\cr
%\noalign{\hrule}
& $q$ &&
$\lambda_{q}(\chi_{4,16})$   &\cr
\noalign{\hrule}
& $2$ &&  $192 \, \alpha$  &\cr
& $3$ && $-497664 \, \alpha +311040 \, \beta +  103935960$ &\cr
& $4$ && $157261824\,\alpha-230584320\,\beta-11214073856$  &\cr
& $5$ &&  $-15895326720\, \alpha +27887846400\, \beta -41100006690$ &\cr
& $7$ &&  $-7428303065088 \, \alpha +12024244661760\, \beta + 44274749992240$ &\cr
& $9$ && $160404429828096\,\alpha-252069882854400\,\beta-377719068915351$  &\cr
} \hrule}
}}
\end{footnotesize}

\smallskip

Harder predicted in \cite{Harder} congruences between elliptic modular forms and 
Siegel modular forms. In the case at hand, the critical value at 
$s=20$ of the $L$-series $\Gamma(s)(2\pi)^{-s} L(f,s)$ of a Hecke 
eigenform $f=\sum_n a(n) q^n$ of weight $34$ for $\SL{2}{\ZZ}$ is 
divisible by the prime $1571$ and we therefore expect for every prime $p$ 
a congruence of the form
$$
\lambda_p(\chi_{4,16})\equiv p^{14}+a(p)+p^{19} \, \bmod \ell\, ,
$$
where $\ell$ is a prime dividing $1571$ in the composite of the fields
${\QQ}(\alpha)$ and ${\QQ}(\sqrt{2356201})$ of eigenvalues $\lambda_p(\chi_{4,16})$
of $\chi_{4,16}$ and $a(p)$ of $f$.
One checks that the norm of  $\lambda_p(\chi_{4,16})-( p^{14}+a(p)+p^{19})$ 
is indeed divisible by $1571$ for $p=2,3,5$ and $7$.

\end{section}

%%%%%%%%%%%%%%%%%%%%%%%%%%%%%%%%%%%%%%%%
\begin{section}{Tables}\label{tables}
In this section we give the Hecke eigenvalues of some modular forms constructed
above using covariants. We checked that these eigenvalues are consistent with the 
results of \cite{BFvdG,BFvdG2,FvdG}. The form $\chi_{8,8}$ generates $S_{8,8}$
and $\chi_{12,6}$ generates $S_{12,6}$. The forms $\chi^{+}_{6,12},
\chi_{6,12}^{-}$ form a basis of the $2$-dimensional space $S_{6,12}$; similar
notation is used
for weights $(10,10)$, $(12,9)$ and $(14,8)$.

\begin{footnotesize}
\smallskip
\vbox{
\bigskip\centerline{\def\quad{\hskip 0.6em\relax}
\def\quod{\hskip 0.5em\relax }
\vbox{\offinterlineskip
\hrule
\halign{&\vrule#&\strut\quod\hfil#\quad\cr
height2pt&\omit&&\omit &&\omit &\cr
%\noalign{\hrule}
& $q$ &&
$\lambda_{q}(\chi_{8,8})$  && $\lambda_{q}(\chi_{12,6})$ &\cr
\noalign{\hrule}
& $2$ &&  $1344$ && $-240$  &\cr
& $3$ && $-6408$ && $68040$  &\cr
& $4$ && $28672$ && $1118464$ & \cr
& $5$ &&  $-30774900$  && $ 1476510$ &\cr
& $7$ && $451366384$ && $-334972400$ &\cr
& $9$ && $-3092097159$ && $-5279708871$ &\cr
& $11$ &&   $13030789224$ && $3580209624$ &\cr
& $13$ && $-328006712228$  && $91151149180$ &\cr
& $17$ && $5520456217764$  && $-11025016477020$ &\cr
& $19$ && $-28220918878760$ &&  $-22060913325080$ &\cr
& $23$ && $79689608755152$  &&  $195863810691120$ &\cr
} \hrule}
}}
\end{footnotesize}

\begin{footnotesize}
\smallskip
\vbox{
\bigskip\centerline{\def\quad{\hskip 0.6em\relax}
\def\quod{\hskip 0.5em\relax }
\vbox{\offinterlineskip
\hrule
\halign{&\vrule#&\strut\quod\hfil#\quad\cr
height2pt&\omit&&\omit &&\omit &\cr
%\noalign{\hrule}
& $q$ && $\lambda_q(\chi_{6,12}^{\pm})$  && $\lambda_{q}(\chi_{10,10}^{\pm})$ &\cr
\noalign{\hrule}
%& && && &\cr
height2pt&\omit&&\omit &&\omit &\cr
& $2$ && $11184 \pm 336\sqrt{601}$ &&  $11400 \pm 120\sqrt{11041}$ &\cr
& $3$ && $-167832 \mp 157248\sqrt {601}$ && $480600 \pm 4320\sqrt {11041}$ &\cr
& $4$ && $-26121728 \pm 5451264\sqrt {601}$ && $149854336 \pm 339840\sqrt {11041}$ &\cr
& $5$ && $-554158500 \pm 77280000\sqrt {601}$ && $-1325429700 \mp 25276800\sqrt {11041}$ &\cr
& $7$ && $-28518281456 \pm 177641856\sqrt {601}$ && $236903926000 \mp 658324800\sqrt {11041}$ &\cr
} \hrule}
}}
\end{footnotesize}

\begin{footnotesize}
\smallskip
\vbox{
\bigskip\centerline{\def\quad{\hskip 0.6em\relax}
\def\quod{\hskip 0.5em\relax }
\vbox{\offinterlineskip
\hrule
\halign{&\vrule#&\strut\quod\hfil#\quad\cr
height2pt&\omit&&\omit &&\omit &\cr
%\noalign{\hrule}
& $q$ &&
$\lambda_q(\chi_{12,9}^{\pm})$  && $\lambda_{q}(\chi_{14,8}^{\pm})$ &\cr
\noalign{\hrule}
height2pt&\omit&&\omit &&\omit &\cr
& 2 && $-6216 \pm 72\sqrt{25249}$ && $-2016\pm96\sqrt{9961}$ &\cr
& 3 && $1074168 \mp 16416\sqrt {25249}$ &&  $2762568\pm18432\sqrt{9961}$ &\cr
& $4$ && $-40492928 \mp 784512\sqrt {25249}$ &&$-80611328\pm1050624\sqrt {9961}$ & \cr
& 5 && $-1795354500 \pm 3600000\sqrt {25249}$ && $-951372900\mp12134400\,\sqrt {9961}$ & \cr
& 7 && $-147859080656 \mp 507187008\sqrt {25249}$ && $87767118544\mp804225024\,\sqrt {9961}$  &\cr
} \hrule}
}}
\end{footnotesize}

\smallskip

Using these tables one can check Harder's conjecture \cite{Harder} 
on the existence of congruences
between Siegel modular forms and elliptic modular forms. In \cite[p. 240]{vdG} 
one finds for the cases of weights $(6,12)$, $(10,10)$, $(12,9)$ and $(14,8)$
predicted congruences with the eigenform $f_{28}^{\pm} \in S_{28}(\Gamma_1)$ modulo $823$,
$157$, $4057$ and $647$. The congruence for the 
eigenvalues for $T_2$ was checked in \cite[p.\ 239]{vdG}
and one can verify that these congruences 
hold for the eigenvalues of $T_3, T_5$ and $T_7$ as well.
For example for $(j,k)=(12,9)$, Harder predicts for every prime $p$ a congruence
$$
\lambda(p)^{\pm} \equiv p^{k-2}+a(p)^{\pm}+p^{j+k-1} \bmod \ell 
$$
for the Hecke eigenvalues of $\chi_{12,9}^{\pm}$ with $a(p)^{\pm}$ the Hecke eigenvalues
of the eigenform $f_{28}^{\pm} \in S_{28}(\Gamma_1)$
with $\ell$ a prime dividing $4057$ in the composite of the fields of eigenvalues
of $f_{28}^{\pm}$ and $\chi_{12,9}^{\pm}$.
One checks that the norm of 
$$
147859080656+507187008\sqrt{25249}+7^7- 87695981800-809077248 \sqrt{18209} + 7^{20}
$$
in ${\QQ}(\sqrt{25249},\sqrt{18209})$ 
is divisible by $4057$.
\end{section}
%%%%%%%%%%%%%%%%%%%%%%%%%%%%%%%%%%%%%%%%

\end{document}